\documentclass[11pt]{amsart}
\usepackage{geometry}              
\usepackage{graphicx}
\usepackage{color}
\usepackage{amssymb}
\usepackage{epstopdf}
\usepackage{amsaddr}
\usepackage{lineno}
 \usepackage[small]{caption}
 \usepackage[section]{placeins}
\newcommand{\black}{\color{black}{}}

\newcommand{\blue}{\color{blue}{}}

\definecolor{orange}{rgb}{0.8, .2, 0.0}

  \newcommand{\BF}{\boldmath}

\title[Many coexisting attractors, almost conservative H\'enon...]{Many coexisting attractors, a case study of the almost-conservative H\'enon map}

\author{Corrado Falcolini}
\address{Department of Architecture, Roma Tre University}
\author{Laura Tedeschini-Lalli}
\address{IMACS (Int'l Association for Mathematics and Computer Simulations)}
\author{James A. Yorke}
\address{University of Maryland College Park}

\email{corrado.falcolini@uniroma3.it, laura.tedeschini.lalli@uniroma3.it,
yorke@umd.edu}
\date{June 8th, 2026}                            

\begin{document}

\begin{abstract}
For dynamical systems in the plane, there can be many periodic attractors coexisting in a bounded region. They become easier to find in systems with small dissipation, which we call ``almost-conservative''.
We ask what happens when there are many periodic attractors. That is the vague question we start with. For a test study, we chose the H\'enon map with a tiny dissipation. We tuned the other parameter to yield a case with 50 attracting periodic orbits. They have a total of 4259 periodic points. We describe how these orbits can be organized into families.
In addition to two low-period orbits, the remaining 48 orbits can be classified into three families, which we describe in detail.
\end{abstract}

\maketitle

\black
\section{Introduction and Notation}

 The phenomenon of multistability is a hallmark of complex nonlinear dynamics, where
systems exhibit a high sensitivity to initial conditions. Research has shown that even
simple planar maps can possess more than 100 coexisting low-period periodic attractors.
Understanding these structures is vital
for the effective control of complex, almost conservative systems. (See \cite{hens2015extreme},\cite{FGHY96}, \cite{FGPY98}) 

 It is well known that dynamical systems can have many periodic attractors or even infinitely many, and that they can be investigated with constructive numerical procedures.
 \cite {yamaguchi2015finding} 
 \cite {TY86}
 
\begin{figure}[!htbp]
\centering
\includegraphics[scale=0.60]{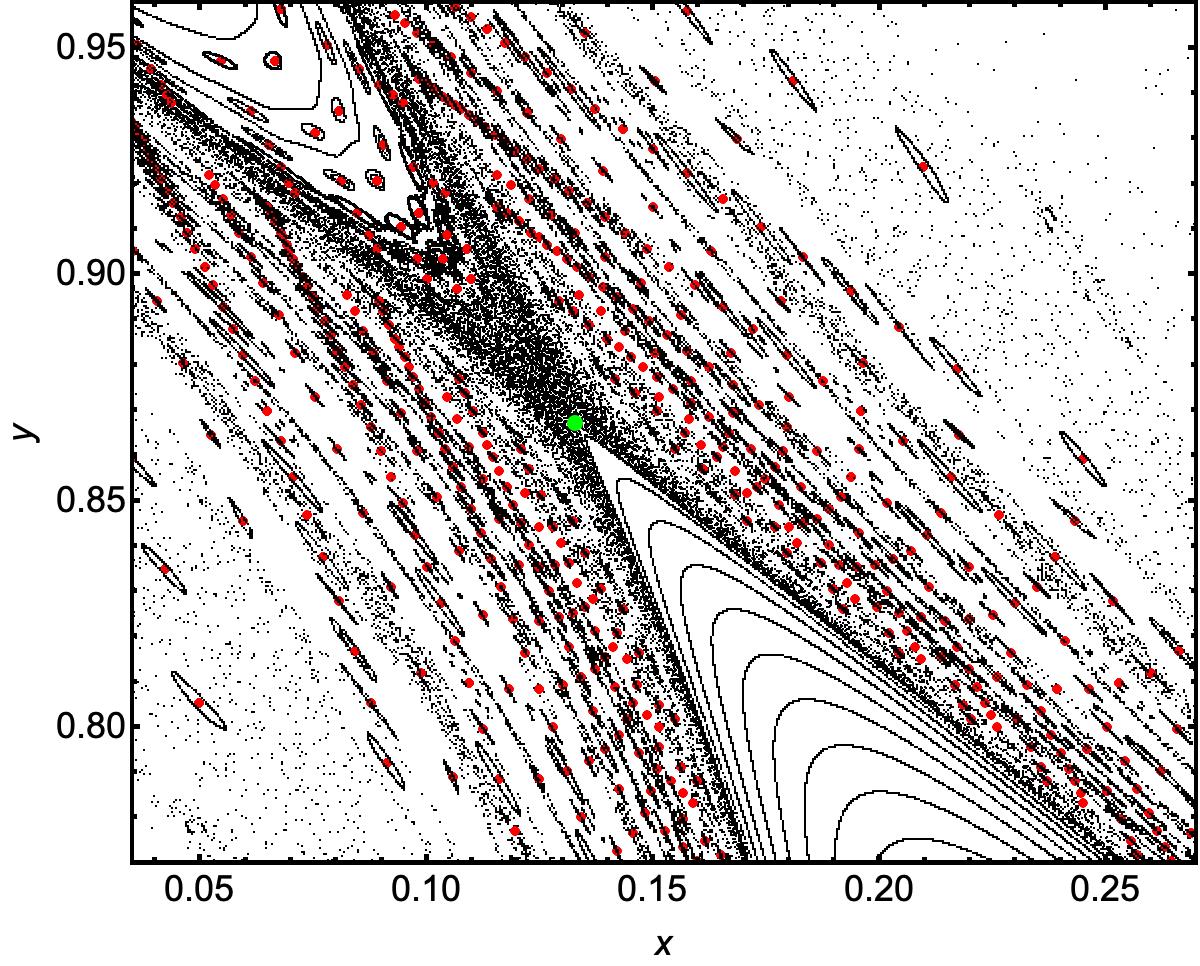}
\caption{{\bf
Elliptic orbits at $\bf b=1$ and $\bf a=1.0176$.} This is the area preserving case. The central green point of this figure is a period-3 saddle point. 
The red dots represent some of the many points that are on elliptical periodic orbits. 
 The black dots are on trajectories of several initial points, each iterated $20{,}000$ times. Some of their orbits trace out closed quasi-periodic curves, and some are in chaotic zones.
}
\label{conservative}
\end{figure}

We study the periodic attractors of the H\'enon map $T_{a,b}(x,y)=(a-x^2-b\ y,x)$,  where $b$ is the Jacobian of the map. For $b=1$, the map is area-preserving and serves as a model for conservative systems (see Figure \ref{conservative}). We investigate the almost conservative case $b=1-10^{-5}$. The dissipation is very small, and the map can have periodic attractors.
For $b>0$, the map is orientation-preserving.
We report on the value of $a=1.0176$, which appears to have a particularly rich assortment of attracting periodic orbits. Here we find three different families of coexisting attractors that have a total of 48 coexisting attracting periodic orbits, see Figure $|2|$. There are also an attracting fixed point and an attracting period-3 orbit.

\vskip0.3cm
\FloatBarrier
\begin{figure}[htbp]
\centering
\includegraphics[scale=0.24]{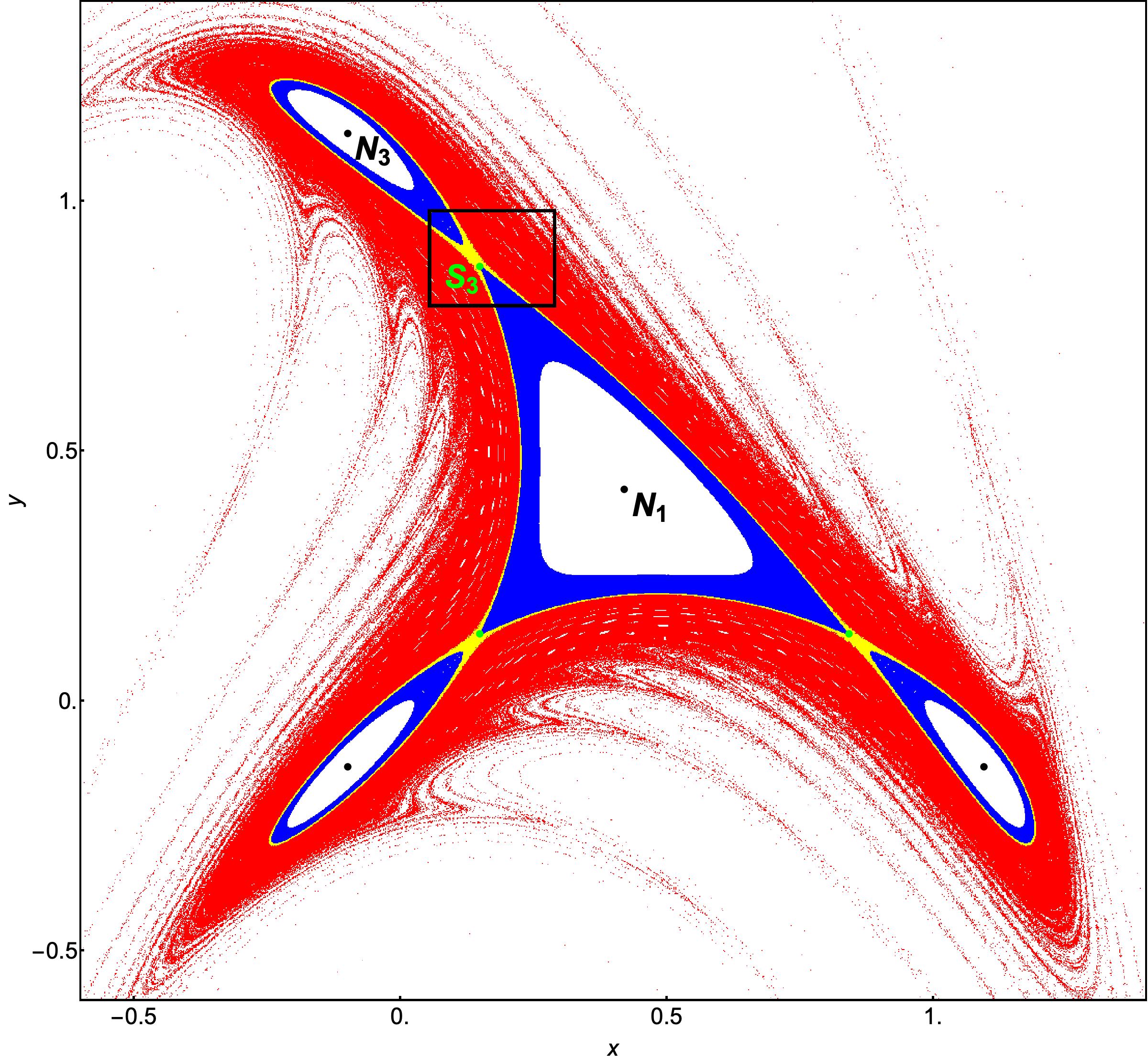} \put (10,25) {(A)} \quad \includegraphics[scale=0.24]{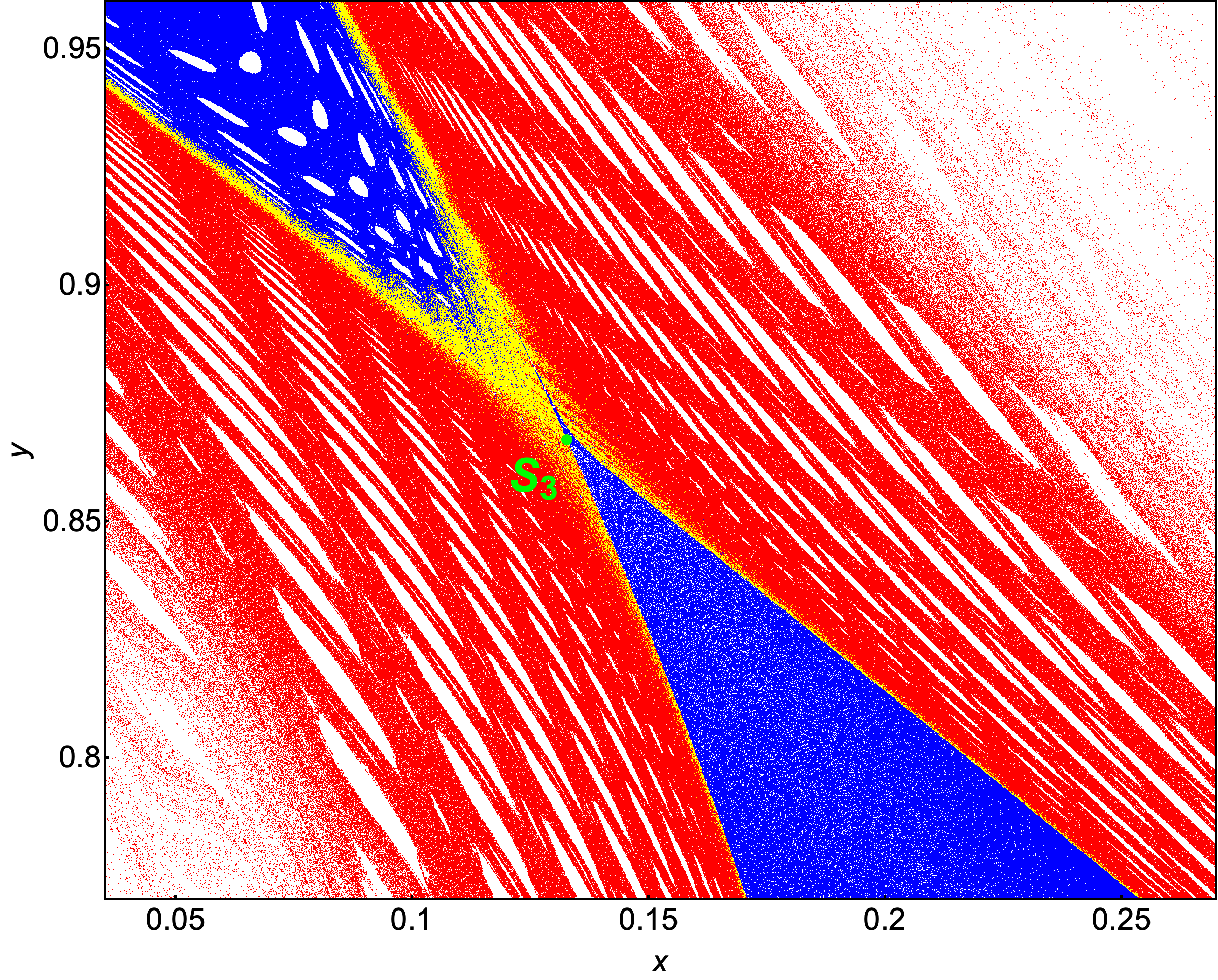} \put (10,25) {(B)} 
\caption{{\bf\BF  Stable and unstable manifolds of a period-3 saddle of the H\'enon map with $\bf a=1.0176$ and $\bf b=1-10^{-5}$.}
$N_1$ is an attracting fixed point. 
$N_3$ is a point on a period-3 attractor. $S_3$ is a point on a period-3 saddle. 
Red points are points on the stable manifold of $S_3$.  Blue points are points on the unstable manifold of $S_3$. 
Yellow pixels contain some of the homoclinic points; these are pixels that contain both stable manifold points and unstable manifolds of $S_3$.
One branch of this stable manifold of $S_3$ is in the basin of the attracting point $N_1$. The other stable branch has a complex behavior with homoclinic intersections, as described in the text. 
{\bf (B)} Enlargement of (A) around the saddle $S_3$.
Each white ‘‘hole" in a blue or red region contains one point of a periodic attractor and the white hole is part of the basin of the orbit of that point.
}
\label{invariant manifolds}
\end{figure}
\FloatBarrier

In almost conservative systems, some general rotation is usually seen (\cite{duarte2000persistent}). In our case, the entire plane rotates approximately $\frac{2\pi}{3}$ counterclockwise. The study in this paper began by investigating areas in the $(a,b)$ parameter plane hinted at in \cite{AS88}. Here we are much closer to the conservative case. In the process we found the coexisting families presented here.

{\bf \BF Notation and Figure \ref{invariant manifolds}.}
Define 
\
\begin{align}\label{Eq:main}
T:=T(x,y)= T_{a,b}(x,y)=(a-x^2-by,x) \text{ where }
a=1.0176,b=1-10^{-5}
\end{align}
The central example of this paper and the example in this figure have those values. 

The basin of attraction of each attracting periodic orbit has infinite area. That follows from the fact that the map shrinks the area $A$ of every region by the factor $b<1$ but the image of a basin is preserved. Therefore, if $A$ is the area of a basin, $A=bA.$ Since a basin has positive area $A>0$, it follows that $A=\infty.$ The pictures we present do not show the infinite area, due to limited resolution and finite computation times.   {These basins include narrow, long filaments whose computation is extremely difficult. The fact that its area is infinite also implies that each basin is unbounded.}

Figure \ref{invariant manifolds}A displays the attracting fixed point $N_1$, a  period-three saddle orbit including the point $S_3$ (green), and a period-three attractor including the black point $N_3$.

Red points are part of the stable manifold of the periodic saddle points of the $S_3$ orbit.
Blue points are part of the unstable manifold of that orbit.
Both are partial plots of those manifolds. If more of the manifolds were plotted, much of the white regions would become red or blue. Each would extend further into its adjacent white areas. We color a pixel yellow when it contains both a point of that stable manifold and a point of that unstable manifold.  The notation is used throughout the paper.
 
 Points in the area of Figure \ref{invariant manifolds} rotate counter-clockwise around $N_1$ by approximately $\frac{2\pi}{3}$. Hence $(x,y)$ will be close to $T^3(x,y)$. In particular, the orbit of $N_3$ rotates counterclockwise around $N_1$.

In Figure \ref{invariant manifolds}, the manifolds are plotted using the method of \cite {kostelich1996plotting}. We can discern three regions: an outer mostly white region, an inner blue region, and a middle red region. Part of the inner blue region appears as a curved triangle containing $N_1$, whose vertices are the green dots of the period three saddle: it consists of one branch of the unstable manifold of point $S_3$ spiraling towards the stable fixed point $N_1$. The other part of the blue region appears as three buttonholes whose vertices are the green dots of the period three saddle: it consists of only one connected curve, i.e., the other branch of the unstable manifold of point $S_3$. Initially, this second branch loops around the attracting periodic orbit $N_3$, but it then goes further, looping around the other two images of $N_3$, in a complicated way. In the middle region, in red, the two branches of the stable manifold of $S_3$ wind closely.

The coexisting attractors described in this paper lie inside the small white holes. The middle region is where the interesting dynamics lies. We will concentrate on basins of attracting periodic orbits in this middle red region. 

\begin{figure}[!htbp]
\centering
\includegraphics[scale=0.25]{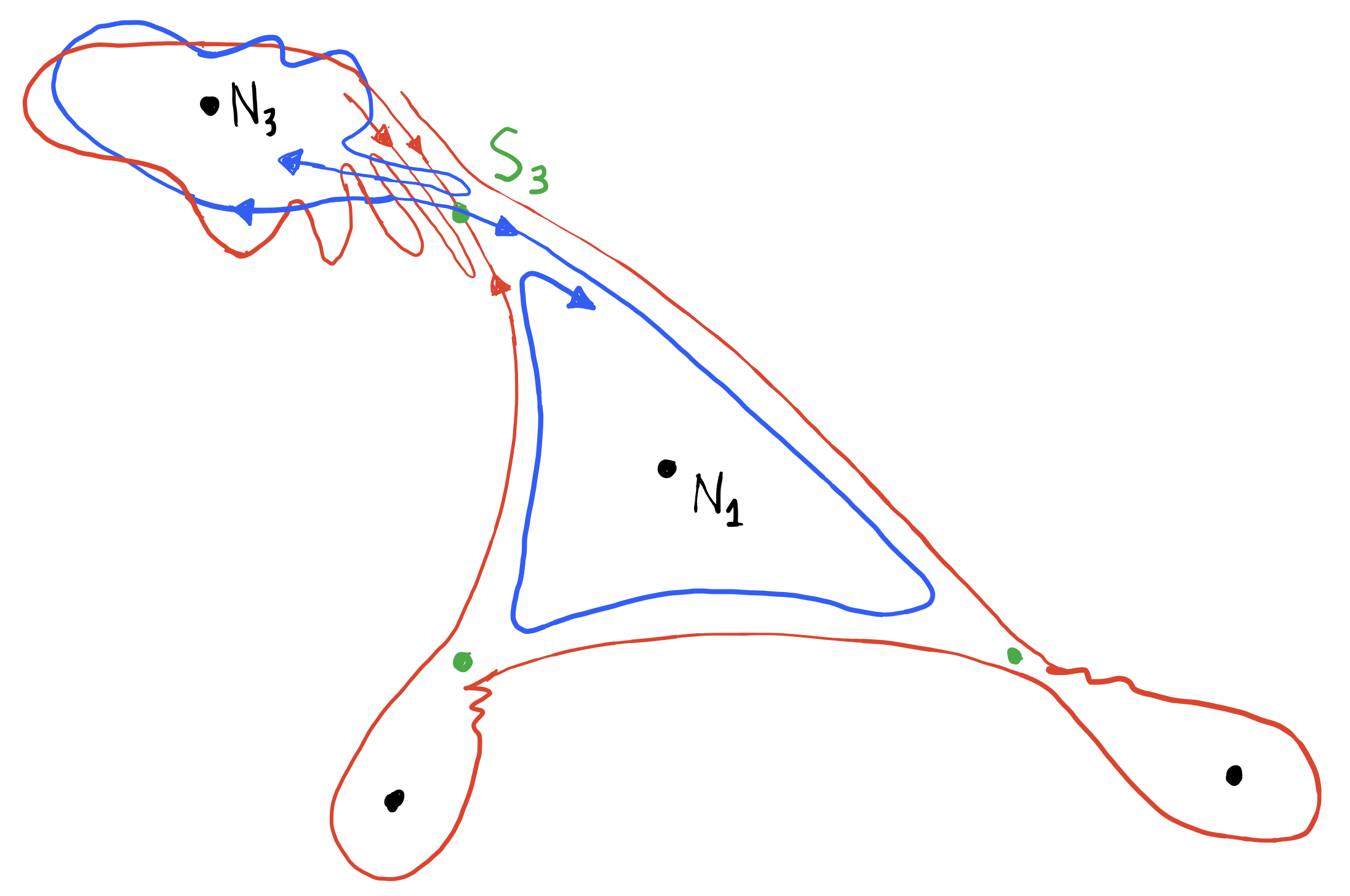}
\caption {{\bf\BF
Qualitative behavior of manifolds of $S_3$.} 
Drawing of stable and unstable manifolds of $S_3$, exaggerating features for visibility.} 
\label{sketchAI}
\end{figure}
\FloatBarrier

As mentioned above, most points in the entire area rotate by approximately $\frac{2\pi}{3}$ around $N_1$. 
In this paper, attracting periodic orbits are organized in families according to their rotation numbers around the point $N_1$

 Some trajectories diverge to $\infty$. 
 The region just described appears to attract all bounded trajectories. We find no other attractors.

\black

\vskip0.5cm

\section {The three families of attracting periodic orbits, their geometric properties in $(x,y)$}

\vskip0.5cm
\noindent {\bf Definition 1}: 

Given a period-$p$ periodic orbit $O_p$ ($O_p=\{Q_0, Q_1, Q_2, \ldots , Q_p\}$, with $Q_0=Q_p$), take vectors $v_k = Q_k - N_1$, from $N_1$ to the points $Q_k$ ($k=0,1,2,\ldots , p$) 
and consider the angles $\alpha (v_k, v_{k+1})\in [0,2\pi)$, between consecutive vectors.
Define its {\bf revolution number ``rev''} (which is an integer) and the (average) {\bf rotation number ``rot''} of $O_p$ as follows.
\begin{align}
rev(O_p)  :=& \frac{1}{2\pi} \sum_{k=0}^{p-1}  \alpha (v_k, v_{k+1});\\
rot(O_p)  :=& \frac{rev(O_p)}{p}; 
\label{1}
\end{align}

\vskip0.5cm

\black
The rotation number of the period-three saddle $S_3$ is $\frac{1}{3}$.  For the periodic attractors we study here, all have a rotation number that is near $\frac{1}{3}$.

Notice that, to preserve geometric information on period $p$ and number of revolutions $rev(O_p)$ of the orbit around $N_1$, the rotation number $rot(O_p)=\frac{rev(O_p)}{p}$ around the 
central point $N_1$ could not be considered in its lowest terms.

\vskip0.5cm

{\bf The three families.}
In this paper we present three families (denoted $F^0, F^1, F^2 $) of coexisting attracting periodic orbits characterized by the form of their rotation number.
For $F^j$, where {\blue $j=0,1,2,$} its period $p$ orbit is denoted by $O_p^j$. It has 
$rev(O_p^j)=(p-j)/3$, which is an integer. See below.
That is, a path that travels along all successive points of a periodic orbit in order, which finally returns to its starting point, will make $(p-j)/3$ 
 revolutions around the fixed point $N_1$ in one period, $p$.

\vskip0.5cm

\noindent {\bf\BF Family $F^0$} has 9 attracting orbits and their periods are $p = 33, 36, \ldots, 57$ (in steps of 3).
Such orbits have the revolution number $rev(O^0_p) = \frac{p}{3}$ and the rotation number $rot(O^0_p)=\frac{1}{3} =\frac{p}{3p}$.

\vskip0.5cm
\noindent 
{\bf\BF Family $F^1$} has 23 attracting orbits and their periods are $p = 37, 40, \ldots, 103$ (in steps of 3). 
Such orbits have revolution number $rev(O^1_p) = \frac{p-1}{3}$ and rotation number $rot(O^1_p) =\frac{p-1}{3p}$. 

\vskip0.3cm
\noindent 
{\bf\BF Family $F^2$} has 16 attracting orbits and their periods are $p = 95, 101, \ldots, 185$ (in steps of 6).
Such orbits have revolution number $rev(O^2_p) = \frac{p-2}{3}$ and rotation number $rot(O^2_p) = \frac{p-2}{3p}$.

\vskip0.3cm
\noindent 
In family $F^j$ for $j=1,2$, the rotation number $rot(O^j_p)$ gets closer to $\frac13$ as $p$ increases. The family $F^j$ has orbits of period $p$ where $p-j$ is divisible by 3. 

\begin{figure}[htbp]
\centering
\includegraphics[scale=0.36]{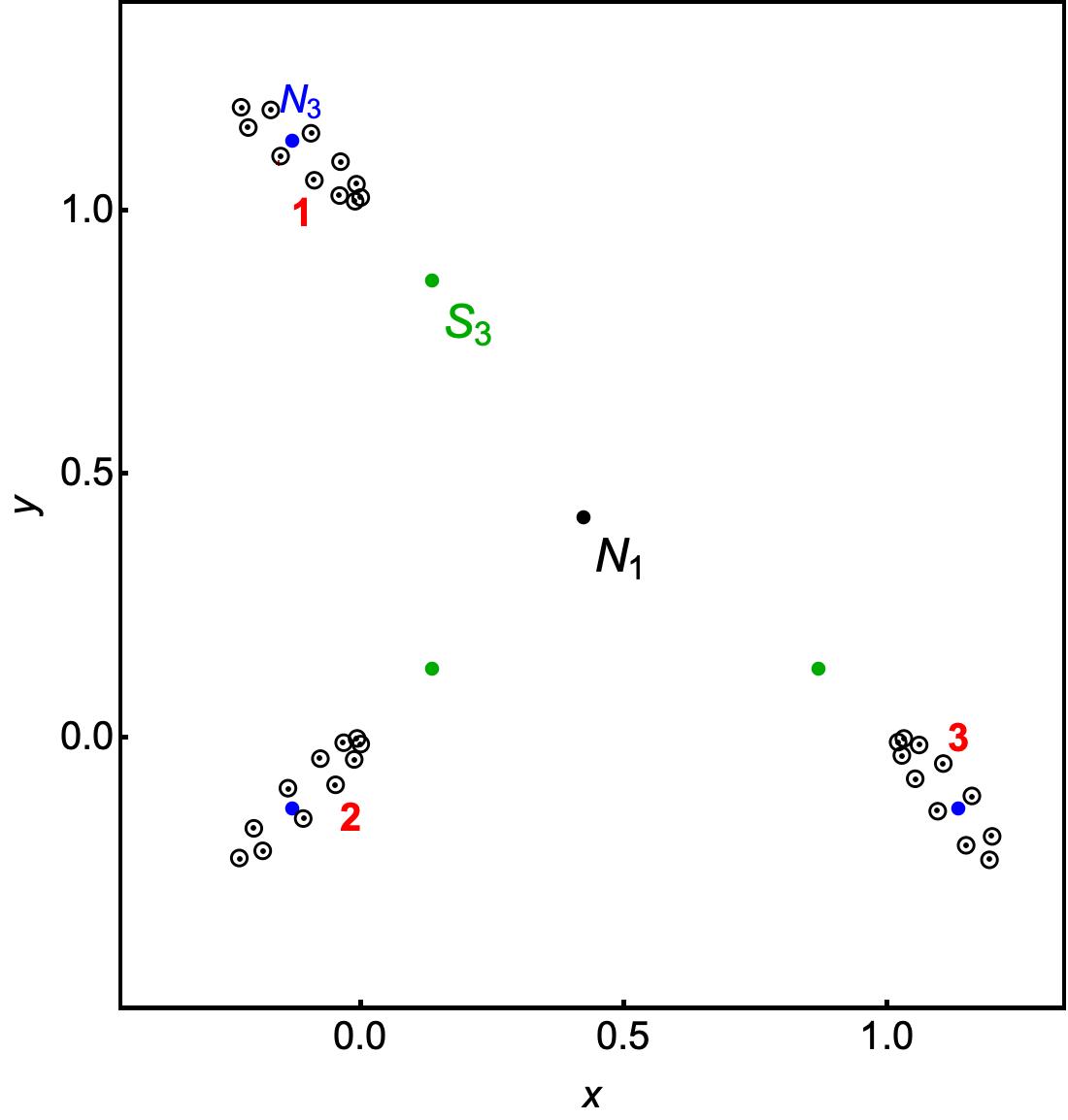} \put (-120,0) {(A)} \quad
\includegraphics[scale=0.36]{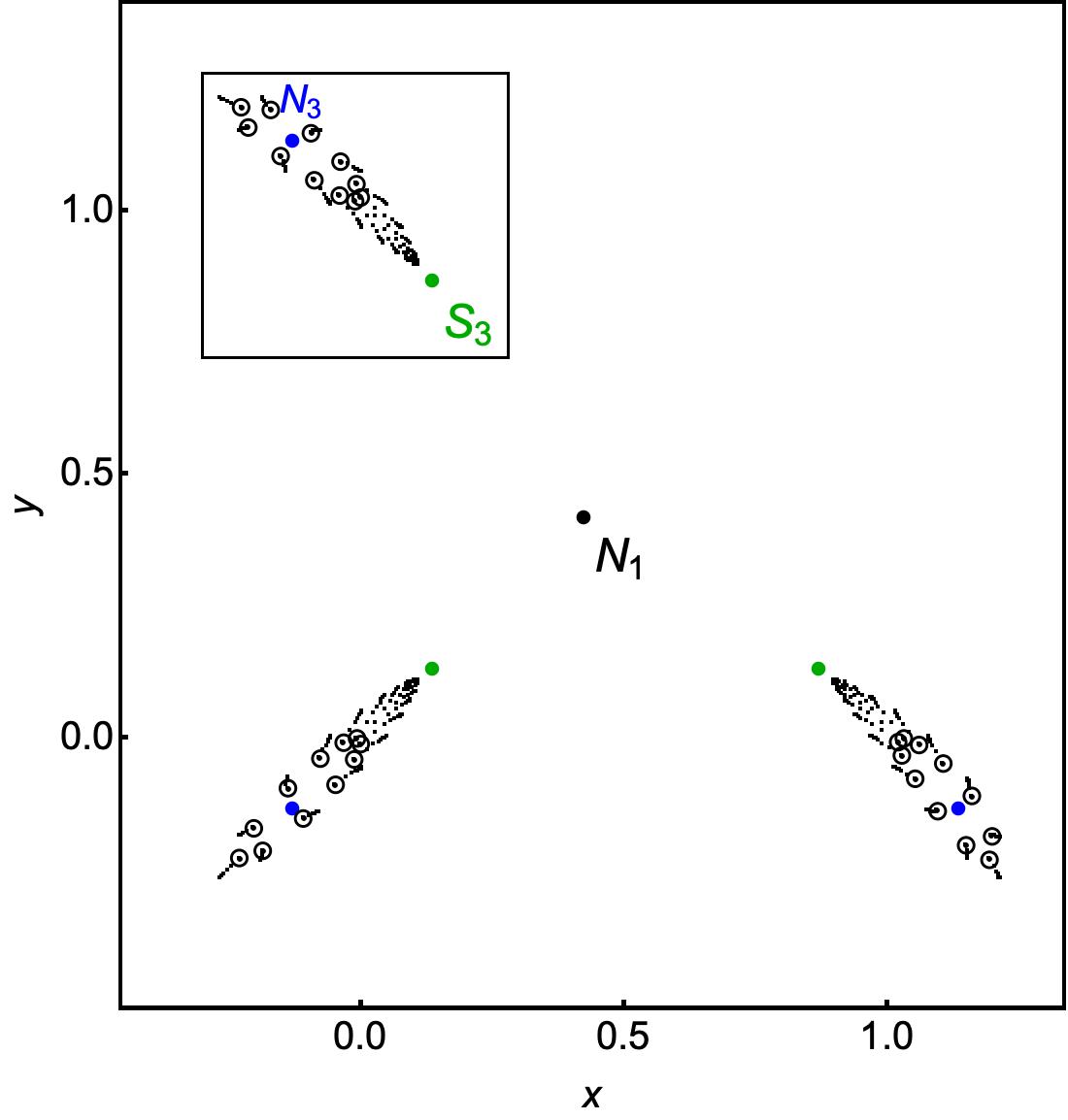} \put (-120,0) {(B)}
\vskip0.3cm
\includegraphics[scale=0.45]{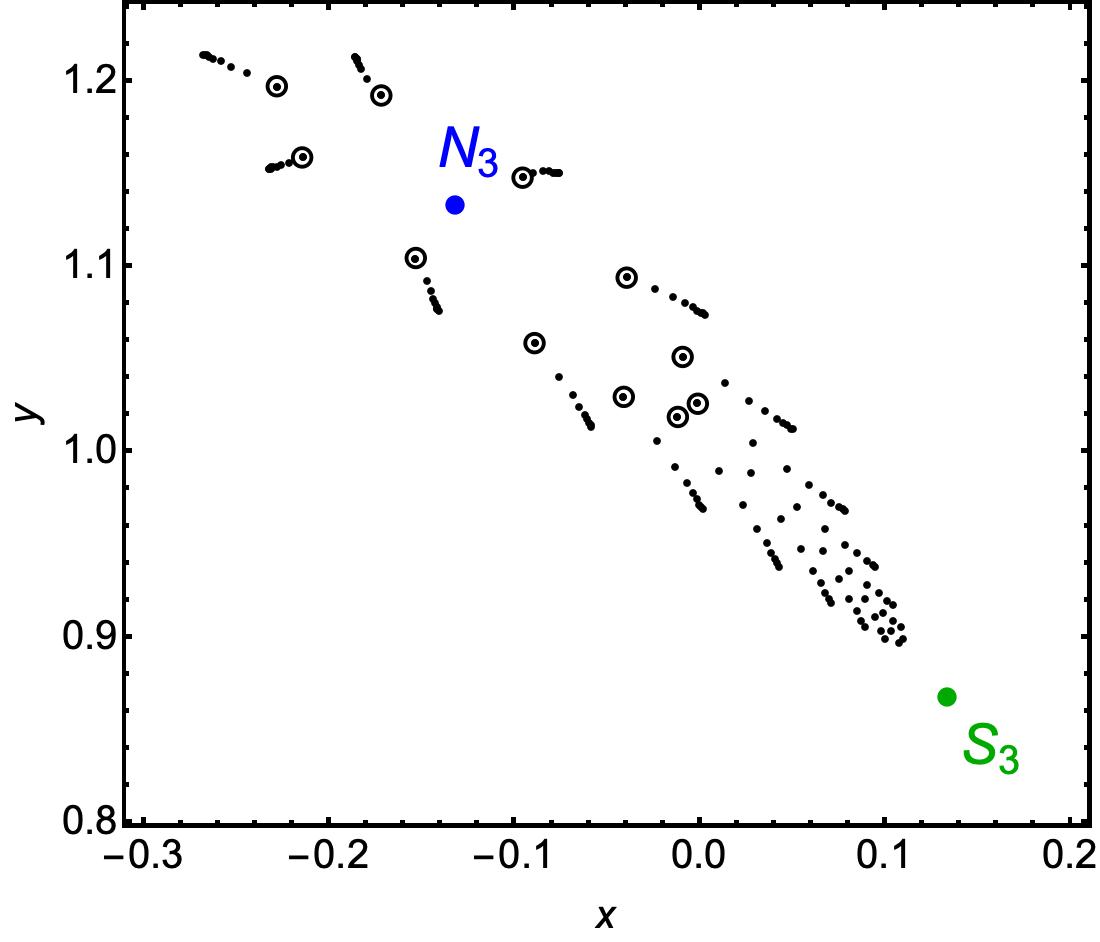} \put (-140,1) {(C)} 

\caption{
{\bf Family $\bf F^0$ consists of 9 periodic sinks with periods 33, 36, 39, \dots, 57, (in steps of 3) with a total of 405 points.} 
$N_1$ (black) is the attracting fixed point, $S_3$ (green) is a period-three saddle point and $N_3$ (blue) is a period-three sink. 
The points of the period-33 orbit are marked by {$\bf \odot$}. 
(A) The period-33 orbit is shown. 
Numbers (red) indicate the order of points on the trajectory, point 1 chosen at random. 
(B) All periodic points of the entire family of $F^0$ sinks (small dots) are shown. (C) It is a blowup from (B).
 }
\label{orbit 33}
\end{figure}

\noindent Notice: the orbits of $F^0$ have rotation number $\frac{1}{3}$.   
\black Each orbit is divisible into three clumps (see Figure \ref{orbit 33}).

Side Note: Each of our families of attracting periodic orbits lies in a larger family of orbits, the others of which are unstable. For example, $F^1$ has attracting orbits of periods $p = 37, 40, \ldots, 103$, and an unstable period-34 orbit. There is also a period-68 attractor, as if the period-34 orbit has period-doubled. We think of this exceptional orbit as part of the family in spirit.

\FloatBarrier

\vskip0.5cm

\noindent For the chosen values of $a=1.0176$ and $b=1-10^{-5}$, the total number of coexisting attracting periodic orbits in the three families is $48$, for a total of  $4255$ points. No other families of attractors were found, other than $N_1$ and the period-three orbit of $N_3$. 

\begin{figure}[htbp]
\centering
\includegraphics[scale=0.36]{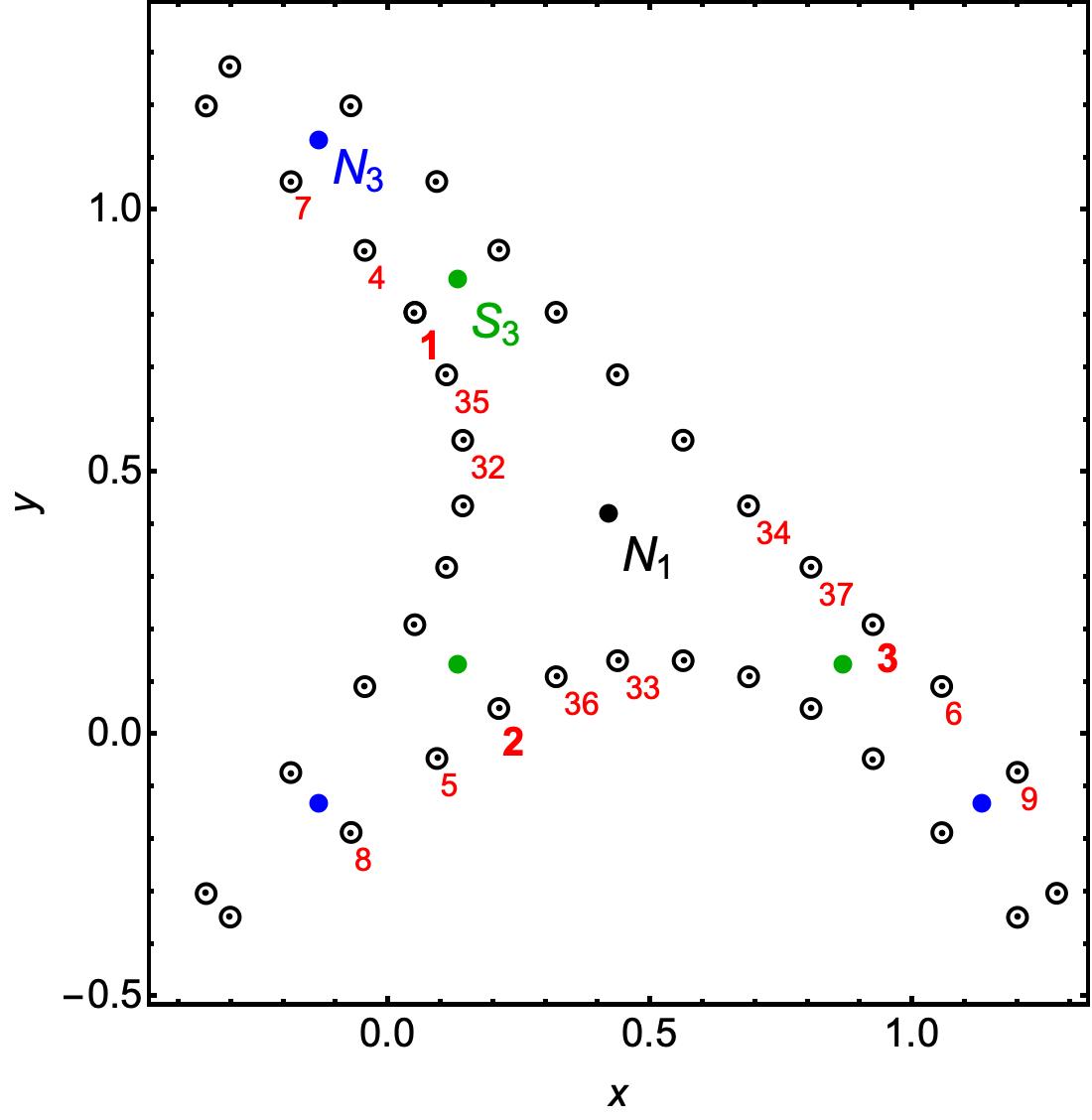} \put (-100,-10) {(A)} \quad
\includegraphics[scale=0.32]{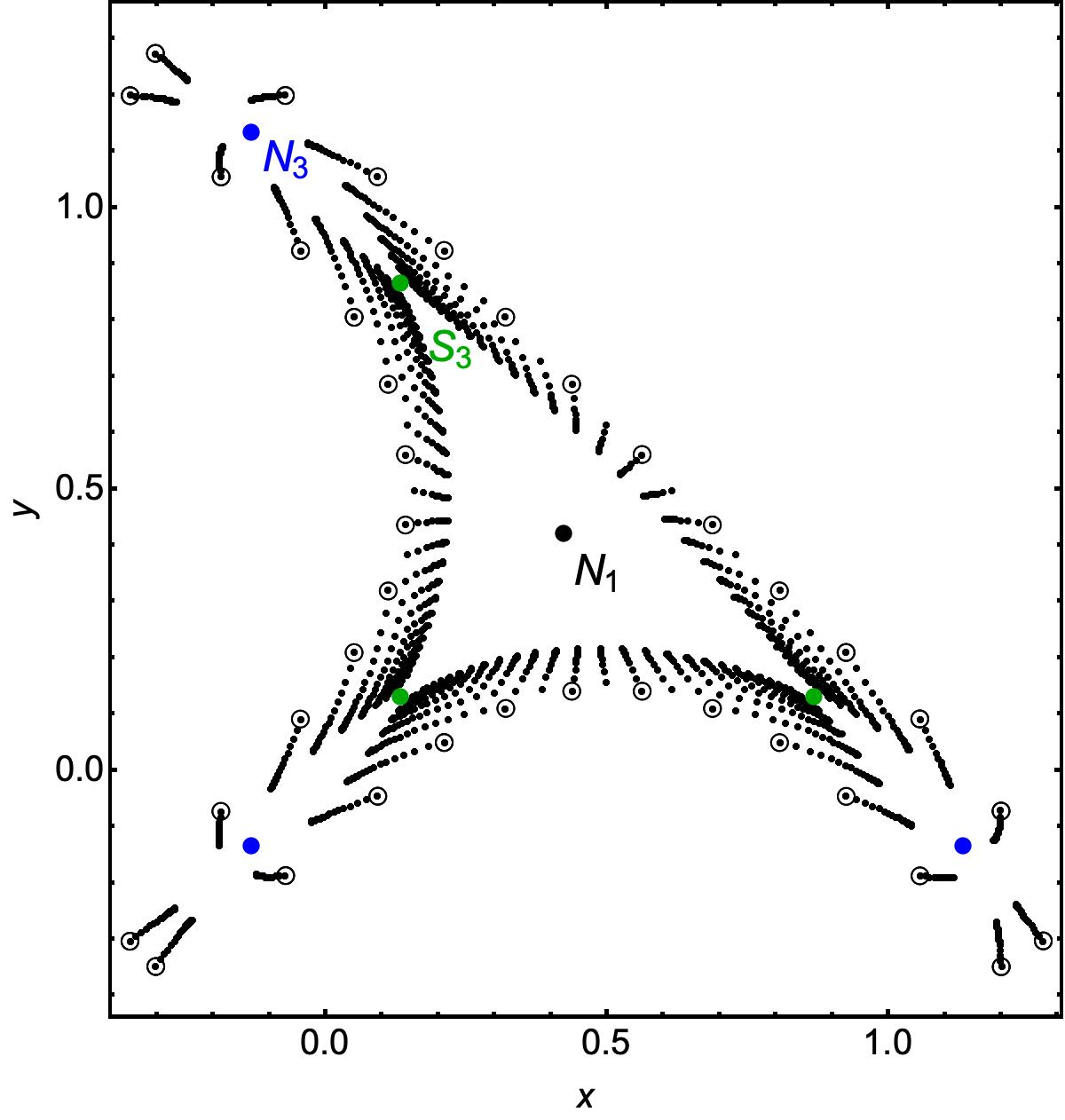} \put (-120,-10) {(B)} 

\caption{
{\bf Family $\bf F^1$ consists of 23 periodic sinks with periods 37, 40, 43, \dots, 103, (in increments of 3) with a total of 1610 points.} 
$N_1$ (black) is the attracting fixed point, $S_3$ (green) is a period-three saddle point and $N_3$ (blue) is a period-three sink. 
The points of the period-37 orbit are marked by {$\bf \odot$}. 
(A) The period-37 orbit is shown. 
Numbers (red) indicate the order of points on the trajectory, with point 1 chosen at random. 
(B) All periodic points of the entire family of $F^1$ sinks (small dots) are shown.
}
\label{orbit 37}
\end{figure}
\FloatBarrier

The lowest period sink of family $F^1$ has period-$37$. 
 The sequence of 37 iterations starting from any of its points makes $12$ revolutions around $N_1$ and has rotation number $rot(O^1_{37})=\frac{12}{37}$ (see Figure \ref{orbit 37}A). All 23 attracting orbits of the $F^1$ family are shown in Figure \ref{orbit 37}. For increasing values of $p$,  their rotation number $rot(O^1_{p})=\frac{p-1}{3p}$ approaches $\frac{1}{3}$, and the periodic orbits seem to approach homoclinic orbits of $S_3$ (see Figure \ref{orbit 37}B).
 
\begin{figure}[!htbp]
\centering
\includegraphics[scale=0.33]{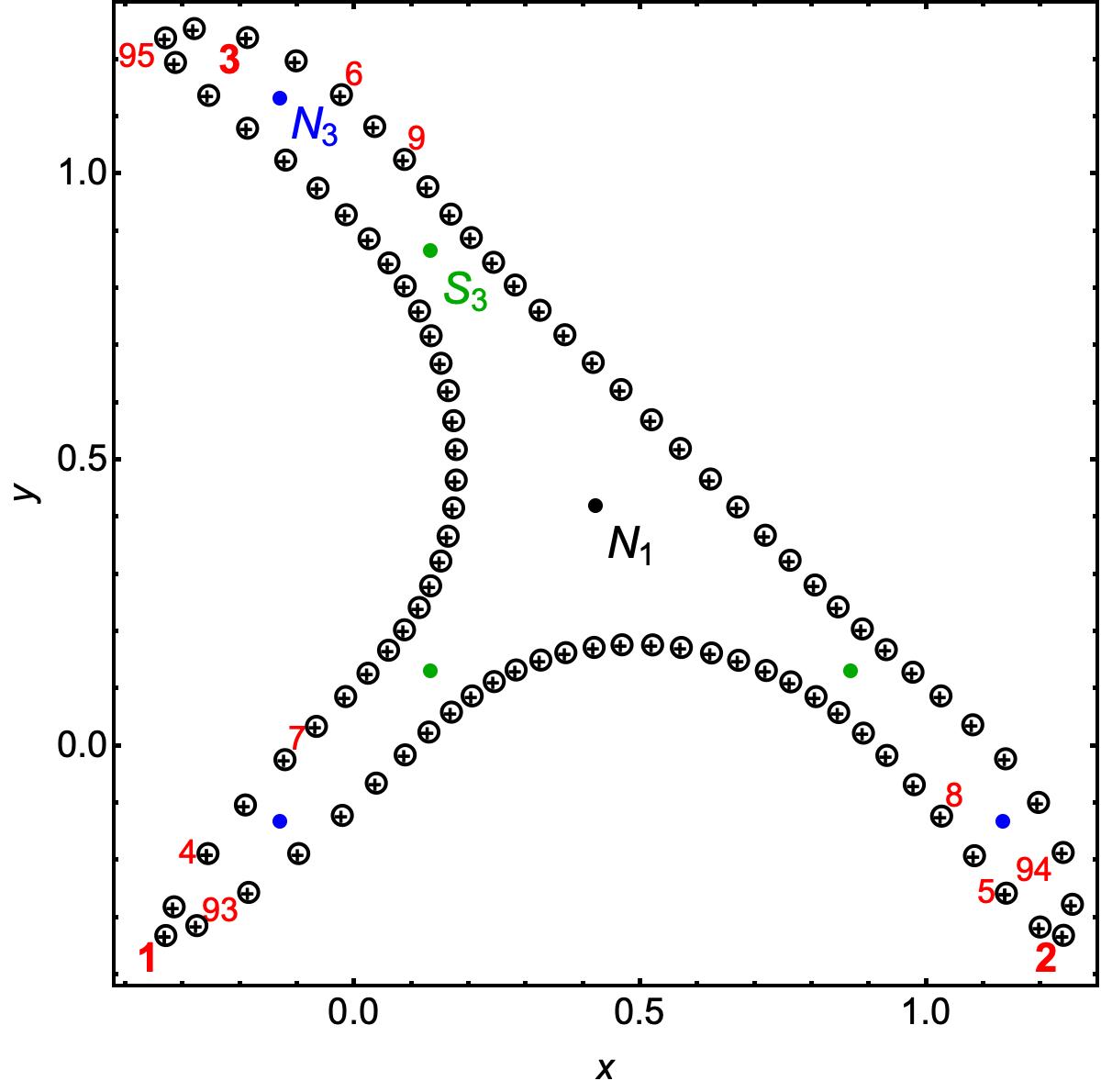} \put (-100,-10) {(A)} \quad
\includegraphics[scale=0.34]{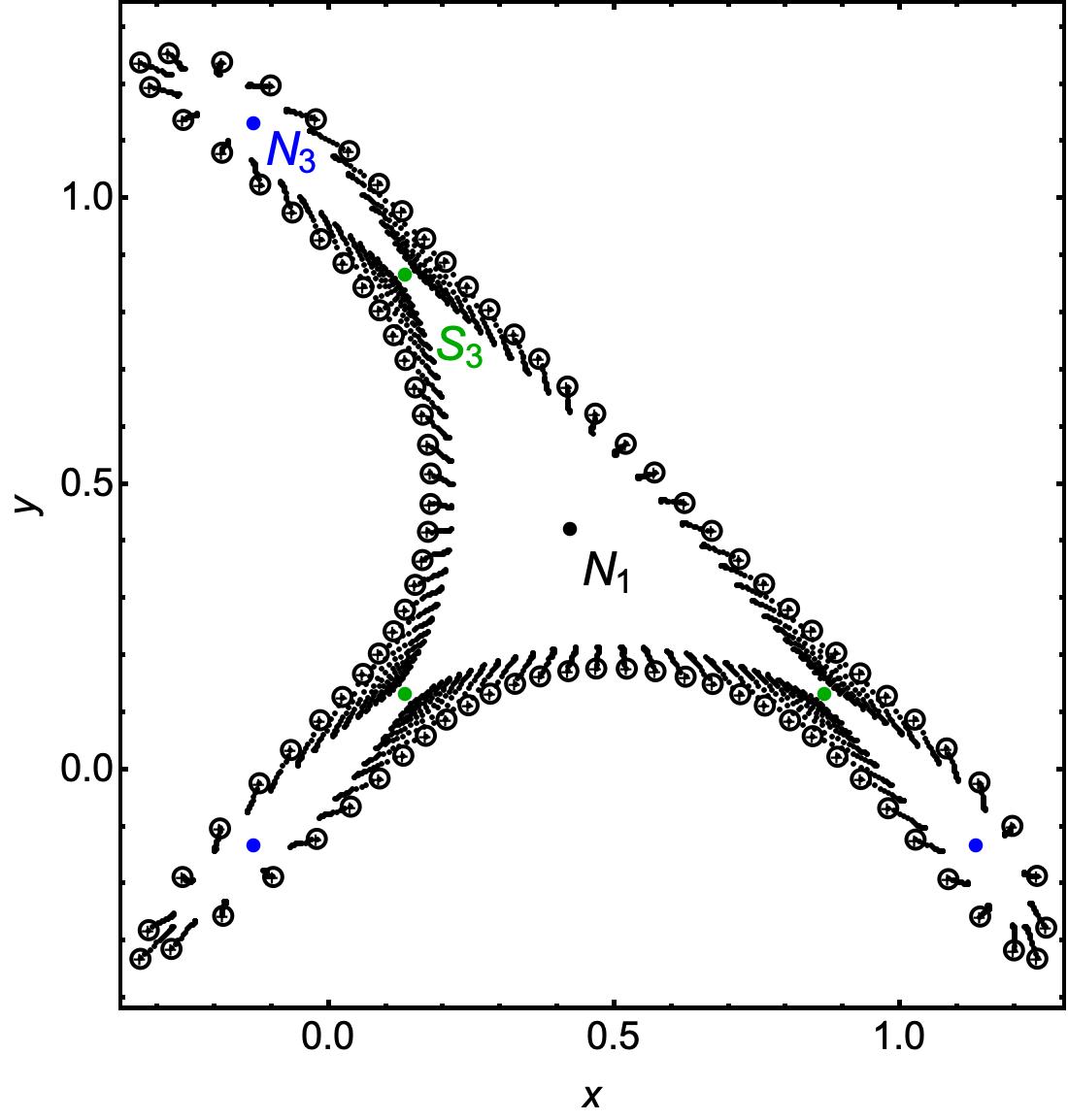} \put (-120,-10) {(B)} 
\caption{
{\bf\BF Family $ F^2$ consists of 16 periodic sinks with periods 95, 101, 107\dots, 185, (in increments of 6) with a total of 2240 points.} 
 $N_1$ (black), $S_3$ (green) and $N_3$ (blue) are as in Figure \ref{orbit 37}. The points of the period-95 orbit are marked with $\bf \oplus$. 
 (A) The period-95 orbit is shown. 
Numbers (red) indicate the order of points on the trajectory, with point 1 chosen at random. 
(B) All periodic points of the entire  $F^2$ family sinks (small dots) are shown. 
}
\label{orbit 95}
\end{figure}
\FloatBarrier

\newpage

 The lowest period sink of family $F^2$ has period-$95$. 
 The sequence of 95 iterations starting from any of its points makes $31$ revolutions around $N_1$ and has rotation number $rot(O^2_{95})=\frac{31}{95}$ (see Figure \ref{orbit 95}A). 
All 16 attracting orbits of the $F^2$ family are shown in Figure \ref{orbit 95}. For increasing values of $p$,  their rotation number $rot(O^2_p) = \frac{p-2}{3p}$ approaches $\frac{1}{3}$, and the periodic orbits seem to approach homoclinic orbits of $S_3$ (see Figure \ref{orbit 95}B).

\begin{figure}[!htbp]
\centering
\includegraphics[scale=0.45]{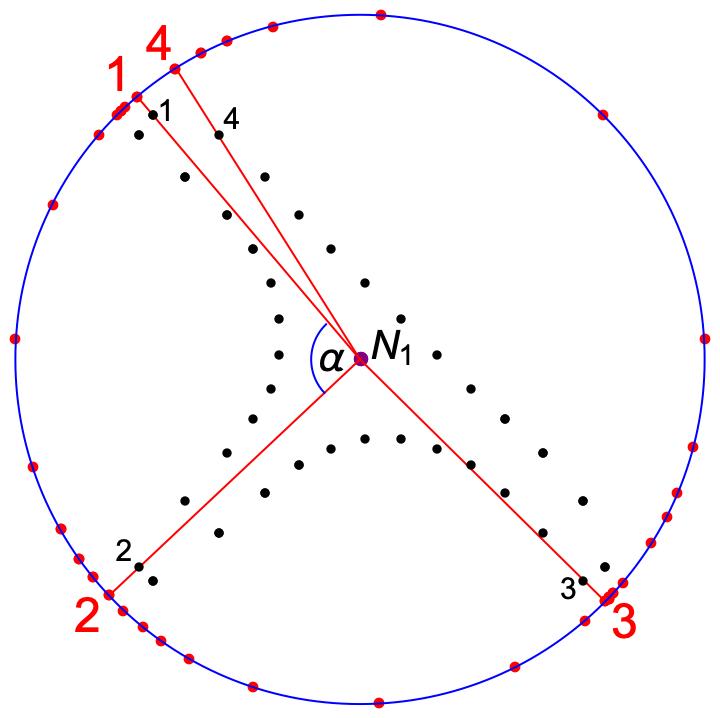} \put (-100,-10) {(A)} \quad \includegraphics[scale=0.45]{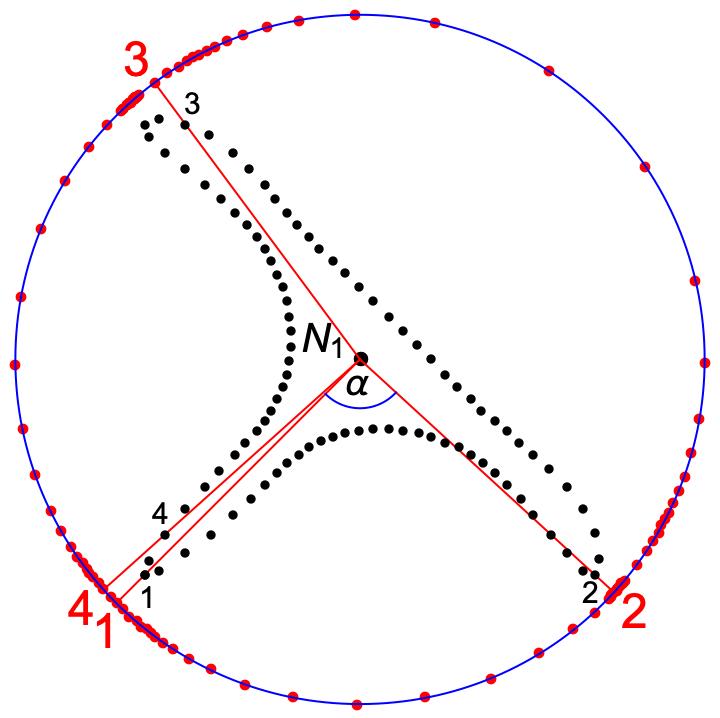} \put (-120,-10) {(B)} 
\caption{{\bf Rotation number of periodic orbits around $\bf N_1$.} 
Consecutive points on each orbit can be viewed as rotating counter-clockwise by an average angle of slightly less than $2\pi/3$, a third of a rotation.
Projection of periodic orbits onto a circle centered at $N_1$. Consecutive points $1, 2, 3, 4$ along the orbit (black) are projected from the center $N_1$ into points $1, 2, 3, 4$ (red) of the circle. (A) The period-$37$ orbit of the $F^1$ family. (B) The period-$95$ orbit of the $F^2$ family.}
\label{cerchio}
\end{figure}
\FloatBarrier

\vskip0.5cm
\noindent In Figure \ref{cerchio}A the period-$37$ orbit is projected on a circle centered at $N_1$. In Figure \ref{cerchio}B the period-$95$ orbit is projected on the same circle centered at $N_1$. Angle $\alpha$ measures the rotation of the map around $N_1$ (see (\ref{1})).  Many of the points of the family are near $S_3$ and its images.

\begin{figure}[!htbp]
\centering
\includegraphics[scale=0.34
]{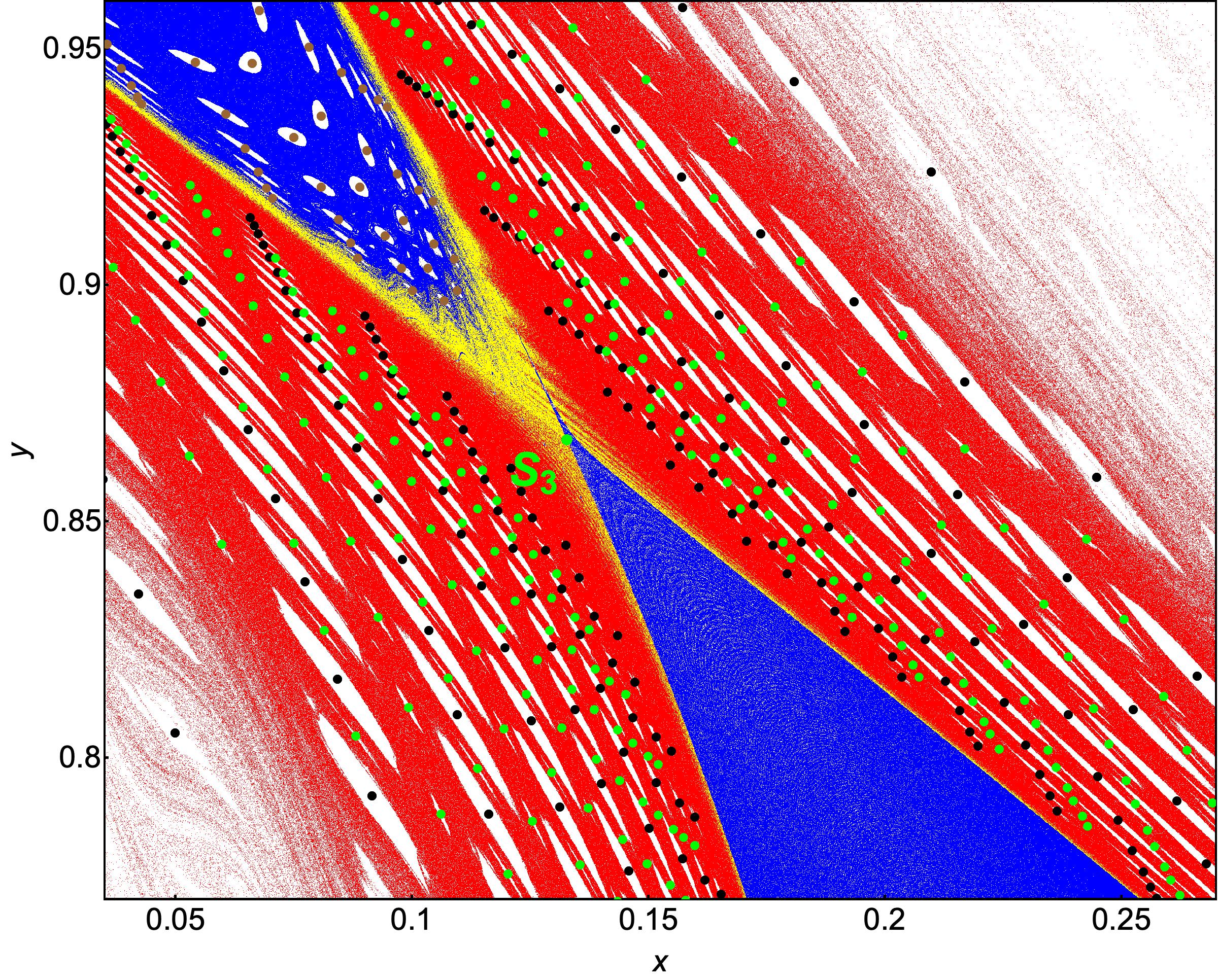}
\caption{{\bf Figure \ref{invariant manifolds}B, shown with all periodic sinks.} Points of the periodic attracting orbits of all three families are shown. They lie in small white oval regions that are part of their basins of attraction. Black dots are periodic attractor points from family $F^1$, green dots from family $F^2$, and brown dots from family $F^0$.
 All 48 orbits of the three families $F^1, F^2, F^0$ are represented by at least one point. The other points lie outside the frame of the picture.}
\label{sovrapposizione}
\end{figure}
\FloatBarrier

\vskip0.5cm
\noindent In Figure \ref{sovrapposizione}, the families $F^1$ and $F^2$ are superimposed on Figure \ref{invariant manifolds}B, where the attracting points are inside the white holes left by the stable (red) invariant manifold of the saddle $S_3$. The orbits of family $F^0$ are inside the white holes left by the unstable (blue) manifold of $S_3$.
\vskip0.5cm

In a chaotic set, most orbits are buried in infinitely many layers of basins of attractors.
There is a special period in which one saddle orbit is in an outermost or innermost position and is called “accessible” from the outside or from the inside.
See  \cite {AY92} for a precise definition of “accessible”.

The lowest period of the family $F_2$ is 95. In addition to that attractor, there is also a period-95 saddle whose 95 points are paired closely with the 95 points of the attractor. This saddle appears to be accessible from within the basin of infinity. It is also the outermost periodic orbit among the three families.
\vskip0.5cm

\section {Maximal coexistence of attractors}

\begin{figure}[!htbp]
\centering
\includegraphics[scale=0.5]{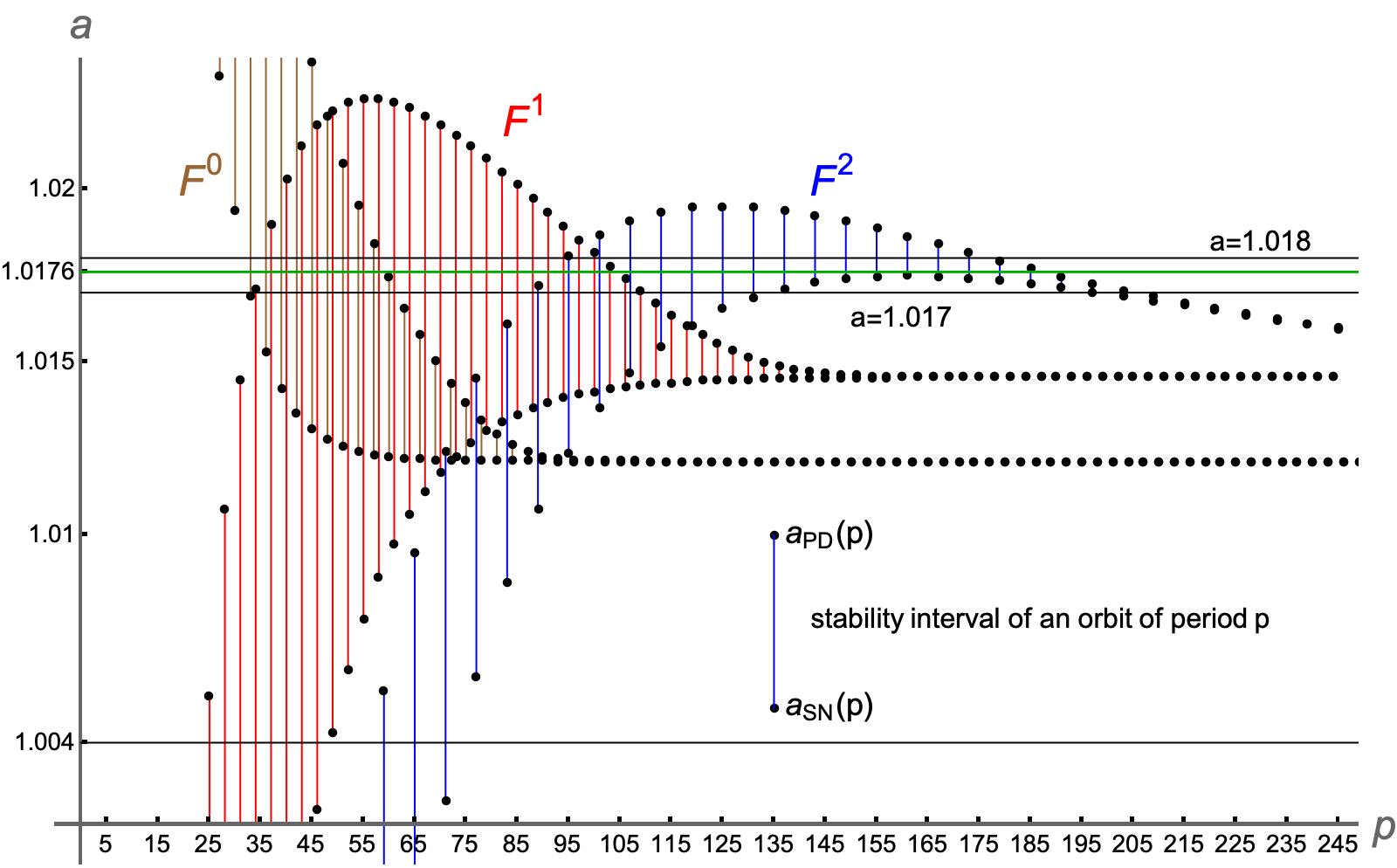}
\caption{{ \bf\BF
Overlapping of stability intervals at $\bf b=1-10^{-5}$ for $a$ near $1.0176$.} 
The horizontal scale shows the period $p$ of a periodic orbit that is stable over a range of $a$ values. 
Each vertical line is a ''stability interval'' for such a periodic orbit. The line is brown for family $F^0$, red for $F^1$ and blue for $F^2$. Each stability interval ends in black dots that represent either a saddle-node bifurcation (bottom) or a period-doubling.
The horizontal line at $a=1.004$ intersects only two of the families, red and blue. 
The value $a=1.0176$ (shown as a horizontal green line) is the subject of this paper. It intersects all three families. Other values of $a$ that intersect all three families have fewer attracting periodic orbits.
The value $a=1.0176$ has been chosen to maximize the total number of coexisting attractors, 48.
For comparison, the horizontal black lines at $a=1.017$ and $a=1.018$ have $43$ and $45$ attractors, respectively.
}
\label{stability intervals}
\end{figure}

We show how the parameter $a$ of the map has been chosen to maximize the number of coexisting attractors across the three families. 
Figure \ref{stability intervals} uses a visualization method introduced in \cite{FT13}.
\vskip0.5cm
 For periodic attractors in this study, as $a$ increases, each orbit arises via a saddle-node bifurcation and then loses stability via a period-doubling bifurcation, consistent with the study in \cite{yorke1983cascades}, as the system is always in the dissipative case. At a given value of $b$, as $a$ increases, the {\it stability interval } of a periodic orbit of period $p$ is therefore defined as the ``$a$'' interval which contains the values of $a$ between the saddle-node and the period-doubling bifurcation points of the orbit. 

\vskip0.5cm

{\bf Figure~\ref{stability intervals} methodology.} The way to find the periodic orbits of high periods discussed in this paper is not obvious. 
Having found three low-period periodic orbits of one of the families, say $F^{1}$, we can find the next orbit as follows. For example, having found the orbits of periods 37,40, and 43. For each, we find the $a$ that is the saddle-node bifurcation point. 
Then choose a point $(x,y)$ from each orbit so that these points are close together. 
From the three $(x,y,a)$ values for each of the three orbits, we predict where the fourth point would be (for the period-46 orbit). We use that as the starting point for a Newton-method search for the saddle-node point $(x,y,a)$. We repeat as needed, finding those points for periods $ 49, 52, \dots, 300$. Figure~\ref{stability intervals} shows these saddle-node points as the bottom points of the red vertical lines. The red lines show the range of $a$ values for which the orbit is stable, ending at the top of the line at a period-doubling point.
We find that this method of rescaling is robust in practice.
See \cite{FT18}. 

\FloatBarrier

The value of $a=1.0176$ (see the green horizontal line in Figure \ref{stability intervals}) has been chosen because each family has several coexisting attractors, $9$ orbits for family $F^0$, $23$ orbits for family $F^1$, and $16$ orbits for family $F^2$. The maximum total number of coexisting attractors could in fact be reached at $a=1.01456508$: $69$ attractors with $11$ orbits of family $F^0$, $54$ orbits of family $F^1$, and $4$ orbits of family $F^2$.

\section {Conclusions}

In this paper, we describe three families of coexisting periodic orbits for the Henon map in the almost conservative case $b=1-10^{-5}$.   In our almost conservative case, the map $T$ rotates the entire plane by approximately $\frac{2\pi}{3}$ counterclockwise around the fixed point $N_1$.

As the periods of each of the three families increase, each family seems to approach points on the homoclinic orbits of the saddle $S_3$, as can perhaps best be seen in Figure~\ref{orbit 33}C.
 We then tuned the parameter $a$ of the map to maximize the number of coexisting periodic attractors. We detected no attractors other than those discussed herein.

\section {Acknowledgments}
 LTL had a long conversation with K.T. Alligood several years ago about the possible implications of the paper \cite{AS88}. This paper stems from curiosity about the general rotational behavior of almost conservative systems.

\bibliographystyle{unsrt}
\bibliography{MAIN.bib}

\end{document}